\documentclass[thmsa,final,11pt]{article}

\usepackage{eurosym}
\usepackage{amsfonts}
\usepackage{amsmath}
\usepackage{amssymb}
\usepackage{mathrsfs} 
\usepackage{enumerate}
\usepackage{graphicx}
\usepackage{subfigure}
\usepackage{color,xcolor}
\usepackage[onehalfspacing,nodisplayskipstretch]{setspace}

\usepackage{xr}

\newtheorem{theorem}{Theorem}

\newtheorem{proposition}{Proposition}

\normalsize
\begin{document}

\title{Fixed-domain asymptotic properties of maximum composite likelihood
estimators for max-stable Brown-Resnick random fields}
\author{Nicolas CHENAVIER\thanks{%
Universit\'{e} du Littoral C\^{o}te d'Opale, 50 rue F. Buisson 62228 Calais.
nicolas.chenavier@univ-littoral.fr}\quad and\quad Christian Y.\ ROBERT%
\thanks{%
1. Laboratory in Finance and Insurance - LFA CREST - Center for Research in
Economics and Statistics, ENSAE, Palaiseau, France; 2. Universit\'{e} de
Lyon, Universit\'{e} Lyon 1, Institut de Science Financi\`{e}re et
d'Assurances, 50 Avenue Tony Garnier, F-69007 Lyon, France. chrobert@ensae.fr%
}}
\maketitle

\begin{abstract}
Likelihood inference for max-stable random fields is in general impossible
because their finite-dimen\-sional probability density functions are unknown
or cannot be computed efficiently. The weighted composite likelihood
approach that utilizes lower dimensional marginal likelihoods (typically
pairs or triples of sites that are not too distant) is rather favored. In
this paper, we consider the family of spatial max-stable Brown-Resnick
random fields associated with isotropic fractional Brownian fields. We
assume that the sites are given by only one realization of a homogeneous
Poisson point process restricted to $\mathbf{C}=(-1/2,1/2]^{2}$ and that the
random field is observed at these sites. As the intensity increases, we
study the asymptotic properties of the composite likelihood estimators of
the scale and Hurst parameters of the fractional Brownian fields using
different weighting strategies: we exclude either pairs that are not edges
of the Delaunay triangulation or triples that are not vertices of triangles.

\strut


\strut

\textit{AMS (2020):} 62G32, 62M30, 60F05, 62H11.
\end{abstract}

\date{}

\section{Introduction}

Gaussian random fields are widely used to model spatial data because their
finite-dimensional distributions are only characterized by the mean and
covariance functions. In general it is assumed that these functions belong
to some parametric models which leads to a parametric estimation problem.
When extreme value phenomena are of interest and meaningful spatial patterns
can be discerned, max-stable random field models are preferred to describe
such phenomena. However, likelihood inference is challenging for such models
because their corresponding finite-dimensional probability density functions
are unknown or cannot be computed efficiently. In this paper we study
composite likelihood estimators in a fixed-domain asymptotic framework for a
widely used class of stationary max-stable random fields: the Brown-Resnick
random fields. As a preliminary, we provide brief reviews of work on maximum
likelihood estimators and on composite likelihood estimators for Gaussian
random fields under fixed-domain asymptotics and present max-stable random
fields with their canonical random tessellations.

\subsection{Maximum likelihood estimators for Gaussian random fields under
fixed-domain asymptotics}

The fixed-domain asymptotic framework is sometimes called infill asymptotics
(Stein (1999), Cressie (1993)) and corresponds to the case where more and
more data are observed in some fixed bounded sampling domain (usually a
region of $\mathbf{R}^{d}$, $d\in \mathbf{N}_{\ast }$). Within this
framework, the maximum likelihood estimators (MLE) of the covariance
parameters of Gaussian random fields have been deeply studied in the last
three decades.

It is noteworthy that two types of covariance parameters have to be
distinguished: microergodic and non-microergodic parameters. A parameter is
said to be microergodic if, for two different values of it, the two
corresponding Gaussian measures are orthogonal (Ibragimov and Rozanov
(1978), Stein (1999)). It is non-microergodic if, even for two different
values of it, the two corresponding Gaussian measures are equivalent.
Non-microergodic parameters cannot be estimated consistently under
fixed-domain asymptotics. No general results are available for the
asymptotic properties of microergodic MLE. Most available results are
specific to particular covariance models.

The initial covariance model that has been studied is the exponential model
with its variance and scale parameters. When $d=1$, only a reparameterized
quantity obtained from the variance and scale parameters is microergodic
(Ying (1991)). It is shown that the MLE of this microergodic parameter is
consistent and asymptotically normal. When $d>1$ and for a separable
exponential covariance function, all the covariance parameters are
microergodic, and the asymptotic normality of the MLE is proved in Ying
(1993). Other results are also given in van der Vaart (1996) and in Abt and
Welch (1998).

The Matern covariance model (Matern (1960)) is very popular in spatial
statistics for its flexibility with respect to the parameterization of
smoothness (in the mean square sense) of the underlying Gaussian field. This
model has three parameters: the variance, the scale and the smoothness
parameters. Zhang (2004) showed that when the smoothness parameter is known
and fixed, not all parameters can be estimated consistently when $d=1,2,3$;
only the ratio of variance and scale parameters (to the power of the
smoothing parameter) is microergodic. Kaufman and Shaby (2013) proved strong
consistency and provided the asymptotic distributions of the microergodic
parameters when estimating jointly the scale and variance parameters (see
also Du et al. (2009) and Wang and Loh (2011) for tapered MLE as well as Loh
et al. (2021) for quadratic variation estimators). For $d=5$, Anderes (2010)
proved the orthogonality of two Gaussian measures with different Matern
covariance functions. In this case, all the parameters are microergodic. The
case $d=4$ is still open.

More recently Bevilacqua et al. (2019) considered the generalized Wendland
(GW) covariance model. They characterized conditions for equivalence of two
Gaussian measures and they established strong consistency and asymptotic
normality of the MLE for the microergodic parameters associated with the GW
covariance model. Bevilacqua and Faouzi (2019) considered the generalized
Cauchy covariance model that is able to separate the characterizations of
the fractal dimension and the long range dependence of the associated
Gaussian random fields. They also characterized conditions for the
equivalence of two Gaussian measures, and established strong consistency and
asymptotic normality of the MLE of the microergodic parameters.

\subsection{Maximum composite likelihood estimators for Gaussian random
fields under fixed-domain asymptotics}

From a theoretical point of view, the maximum likelihood method is the best
approach for estimating the covariance parameters of a Gaussian random
field. Nevertheless, the evaluation of the likelihood function under the
Gaussian assumption requires a computational burden of order $O(n^{3})$ for $%
n$ observations (because of the inversion of the $n\times n$ covariance
matrix), making this method computationally impractical for large datasets.
The composite likelihood (CL) methods rather use objective functions based
on the likelihood of lower dimensional marginal or conditional events (Varin
et al. (2011)). These methods are generally appealing when dealing with
large data sets or when it is difficult to specify the full likelihood, and
provide estimation methods with a good balance between computational
complexity and statistical efficiency.

There is not a lot of results under fixed domain asymptotics for maximum CL
estimators (MCLE). However, Bachoc et al. (2019) studied the problem of
estimating the covariance parameters of a Gaussian process ($d=1$) with
exponential covariance function. They showed that the weighted pairwise
maximum likelihood estimator of the microergodic parameters can be
consistent, but also inconsistent, according to the objective function; e.g.
the weighted pairwise conditional maximum likelihood estimator is always
consistent (and also asymptotically Gaussian). Bachoc and Lagnoux (2020)
considered a Gaussian process ($d=1$) whose covariance function is
parametrized by variance, scale and smoothness parameters. They focused on
CL objective functions based on the conditional log likelihood of the
observations given the $K$ (resp. $L$) observations corresponding to the
left (resp. right) nearest neighbor observation points. They examined the
case where only the variance parameter is unknown and the case where the
variance and the spatial scale are jointly estimated. In the first case they
proved that for small values of the smoothness parameter, the composite
likelihood estimator converges at a sub-optimal rate and they showed that
the asymptotic distribution is not Gaussian. For large values of the
smoothness parameter, they proved that the estimator converges at the
optimal rate.

\subsection{Fixed-domain asymptotics for non-Gaussian random fields}

To the best of our knowledge, there is a few papers that study MLE\ or MCLE
for non-Gaussian random fields under fixed-domain asymptotics. For example,
Li (2013) proposed approximate maximum-likelihood estimation for diffusion
processes ($d=1$) and provided closed-form asymptotic expansion for
transition density. But diffusion processes may not be generalized for $%
d\geq 2$.

Other papers rather considered variogram-based or power variation-based
estimators. Chan and Wood (2004) considered a random field of the form $%
g\left( X\right) $, where $g:\mathbf{R}\mathbb{\rightarrow }\mathbf{R}$ is
an unknown smooth function and $X$ is a real-valued stationary Gaussian
field on $\mathbf{R}^{d}$ ($d=1$ or $2$) whose covariance function obeys a
power law at the origin. The authors addressed the question of the
asymptotic properties of variogram-based estimators when $g\left( X\right) $
is observed instead of $X$ under a fixed-domain framework. They established
that the asymptotic distribution theory for nonaffine $g$ is somewhat richer
than in the Gaussian case (i.e. when $g$ is an affine transformation).
Although the variogram-based estimators are not MLE\ or MCLE, this study
shows that their asymptotic properties can differ significantly from the
Gaussian random field case. Robert (2020) considered a particular class of
max-stable processes ($d=1$), the class of simple Brown-Resnick max-stable
processes whose spectral processes are continuous exponential martingales.
He developed the asymptotic theory for the realized power variations of
these max-stable processes, that is, sums of powers of absolute increments.
He considered a fixed-domain asymptotic setting and obtained a biased
central limit theorem whose bias depends on the local times of the
differences between the logarithms of the underlying spectral processes.

\subsection{Max-stable random fields}

Max-stable random fields appear as the only possible non-degenerate limits
for normalized pointwise maxima of independent and identically distributed
(i.i.d.) random fields with continuous sample paths (see e.g. de Haan and
Ferreira (2006)). The one-dimensional marginal distributions of max-stable
fields belong to the parametric class of Generalized Extreme Value
distributions. Since we are interested in the estimation of parameters
characterizing the dependence structure, we restrict our attention to
max-stable random fields $\eta =(\eta (x))_{x\in \mathcal{X}}$ on $\mathcal{%
X\subset }\mathbf{R}^{d}$ with standard unit Fr\'{e}chet margins, that is,
satisfying%
\begin{equation*}
\mathbb{P}\left[ \eta (x)\leq z\right] =\exp \left( -z^{-1}\right) ,\qquad 
\text{for all }x\in \mathcal{X}\text{ and }z>0\text{.}
\end{equation*}%
The max-stability property has then the simple form 
\begin{equation*}
n^{-1}\bigvee_{i=1}^{n}\eta _{i}\overset{d}{=}\eta
\end{equation*}%
where $(\eta _{i})_{1\leq i\leq n}$ are i.i.d. copies of $\eta $, $\bigvee $
is the pointwise maximum, and $\overset{d}{=}$ denotes the equality of
finite-dimensional distributions. Max-stable random fields are characterized
by their spectral representation (see e.g., de Haan (1984), Gin\'{e} et al.
(1990)): any stochastically continuous max-stable process $\eta $ can be
written as%
\begin{equation}
\eta (x)=\bigvee_{i\geq 1}U_{i}Y_{i}(x),\qquad x\in \mathcal{X}\text{,}
\label{Eq_Spectral_representation}
\end{equation}%
where $(U_{i})_{i\geq 1}$ is the decreasing enumeration of the points of a
Poisson point process on $(0,+\infty )$ with intensity measure $u^{-2}%
\mathrm{d}u$, $(Y_{i})_{i\geq 1}$ are i.i.d. copies of a non-negative
stochastic random field $Y$ on $\mathcal{X}$ such that $\mathbb{E}[Y(x)]=1$
for all $x\in \mathcal{X}$, the sequences $(U_{i})_{i\geq 1}$ and $%
(Y_{i})_{i\geq 1}$ are independent.

The spectral representation $\left( \ref{Eq_Spectral_representation}\right) $
makes it possible to construct a canonical tessellation of $\mathcal{X}$ as
in Dombry and Kabluchko (2018). We define the cell associated with each
index $i\geq 1$ by $C_{i}=\{x\in \mathcal{X}:U_{i}Y_{i}(x)=\eta (x)\}$. It
is a (possibly empty) random closed subset of $\mathcal{X}$ and each point $%
x\in \mathcal{X}$ belongs almost surely (a.s.) to a unique cell (the point
process $\{U_{i}Y_{i}(x)\}_{i\geq 1}$ is a Poisson point process with
intensity $u^{-2}\mathrm{d}u$ so that the maximum $\eta (x)$ is almost
surely attained for a unique $i$). It is noteworthy that the terms \textit{%
cell} and \textit{tessellation} are meant in a broader sense than in
Stochastic Geometry where they originated. Here, a cell is a general (not
necessarily convex or connected) random closed set and a tessellation is a
random covering of $\mathcal{X}$ by closed sets with pairwise disjoint
interiors.

Likelihood inference is challenging for max-stable random fields because
their finite-dimensional probability density functions are unknown or cannot
be computed efficiently. Padoan et al. (2010) proposed to use a
composite-likelihood approach but only discussed asymptotic properties of
the estimators when the data-sites are fixed and when there is a large
number of i.i.d. data replications.

\subsection{Contributions of the paper}

In this paper, we consider the class of spatial max-stable Brown-Resnick
random fields ($d=2$) associated with isotropic fractional Brownian random
fields as defined in Kabluchko et al. (2009). We assume a Poisson stochastic
spatial sampling scheme and use the Poisson-Delaunay triangulation to select
the pairs and triples of sites with their associated marginal distributions
that will be integrated into the CL objective functions (we exclude pairs
that are not edges of the Delaunay triangulation or triples that are not
vertices of triangles of this triangulation). Note that using the Delaunay
triangulation is relatively natural here since we only use the distributions
of pairs and triples. Moreover, the Delaunay triangulation appears to be the
most \textquotedblleft regular\textquotedblright\ triangulation in the sense
that it is the one that maximises the minimum of the angles of the triangles.

We study for the first time the asymptotic properties of the MCLE of the
scale and Hurst parameters of the max-stable Brown-Resnick random fields
under fixed domain asymptotics (for only one realization of a Poisson point
process). Pairwise and triplewise CL objective functions (considering all
pairs and triples) have been proposed for inference for max-stable
processes, but the properties of the MCLE have only been studied when the
sites are fixed and when there is a large number of independent observations
over time of the max-stable random field (see, e.g., Blanchet and Davison
(2011), Davison et al. (2012) or Huser and Davison (2013)). Note that the
tapered CL estimators for max-stable random field excluding pairs that are at a
too large distance apart have also been studied in Sang and Genton (2014)
(here again with independent observations), but this is the first time that
a Delaunay triangulation is used to select the pairs and triples.

To obtain the asymptotic distributions of the MCLEs, we proceed in several
steps. First we consider sums of square increments of an isotropic
fractional Brownian field on the edges of the Delaunay triangles and provide
asymptotic results using Malliavin calculus (see Theorem \ref%
{Prop_Conv_V_n_Gaussian}). Zhu and Stein (2002) also studied sums of
generalized variations for this random field but assumed data-sites on a
regular grid. Second we consider sums of square increments of the pointwise
maximum of two independent isotropic fractional Brownian fields and show
that the asymptotic behaviors of the sums now depend on the local time at
the level $0$ of the difference between the two fractional Brownian fields
(see Theorem \ref{cor:twotrajectories}). Third we generalize these results
to the max-stable Brown-Resnick random field which is built as the pointwise
maximum of an infinite number of isotropic fractional Brownian fields (see
Theorem \ref{prop:BRtrajectories}). Using approximations of the pairwise and
triplewise CL objective functions, we derive the asymptotic properties of
the MCLEs (see Theorem \ref{Prop:Asym_Prop_CL_Est}).

\strut

The family of stationary Brown-Resnick random fields defined in Kabluchko et
al. (2019) is presented in Section 2. We also provide the asymptotic
distributions of pairs and triples as the distances between sites tend to
zero. In Section 3, we introduce the randomized sampling scheme and define
the CL estimators of the scale and Hurst parameters. Our main results are
stated in Section 4. The proofs and some intermediate results are deferred
into a Supplementary Material.

\section{The max-stable Brown-Resnick random fields}

\subsection{Definition of the max-stable Brown-Resnick random fields}

This paper concerns the class of max-stable random fields known as
Brown-Resnick random fields. This class of random fields is based on
Gaussian random fields with stationary increments and was introduced in
Kabluchko et al. (2009). Recall that a random process $\left( W\left(
x\right) \right) _{x\in \mathbf{R}^{d}}$ is said to have stationary
increments if the law of $\left( W\left( x+x_{0}\right) -W\left(
x_{0}\right) \right) _{x\in \mathbf{R}^{d}}$ does not depend on the choice
of $x_{0}\in \mathbf{R}^{d}$. A prominent example is the isotropic
fractional Brownian field where $W\left( 0\right) =0$ a.s. and
semi-variogram given by $\gamma \left( x\right) =\text{var}\left( W\left( x\right)
\right) /2=\sigma ^{2}\left\Vert x\right\Vert ^{\alpha }/2$ for some $\alpha
\in (0,2)$ and $\sigma ^{2}>0$, where $\left\Vert x\right\Vert $ is the
Euclidean norm of $x$. The parameter $\sigma $ is called the scale parameter
while $\alpha $ is called the range parameter ($H=\alpha /2$ is also known
as the Hurst parameter and relates to the H\"{o}lder continuity exponent of $%
W$). It is noteworthy that $W$ is a self-similar random field with linear
stationary increments as presented in Definition 3.3.1 of Cohen and Istas
(2013) and it differs from the fractional Brownian sheet which is a
self-similar random field with stationary rectangular increments (see e.g.
Section 3.3.2 of the same book). Functional limit theorems for generalized
variations of this fractional Brownian sheet have been studied in Pakkanen
and Reveillac (2016), but these theorems cannot be extended to the isotropic
fractional Brownian field whose rectangular increments are not stationary.

In this paper we consider spatial max-stable random fields ($d=2$) and
assume that the random field $Y$ introduced in the spectral representation $%
\left( \ref{Eq_Spectral_representation}\right) $ has the following form%
\begin{equation*}
Y\left( x\right) =\exp \left( W\left( x\right) -\gamma \left( x\right)
\right) ,\qquad x\in \mathbf{R}^{2}.  \label{Eq_Def_Y}
\end{equation*}%
With this choice, $\eta $ is a stationary random field while $W$ is not
stationary but has (linear) stationary increments (see Kabluchko et al.
(2009)).

\subsection{Pairwise joint distributions and asymptotic score contributions}

Let us consider two sites $x_{1},x_{2}\in \mathbf{R}^{2}$ and denote by $%
d=\left\Vert x_{2}-x_{1}\right\Vert $ the distance between these sites. Let $%
z_{1},z_{2}\in \mathbf{R}_{+}$, $a=\sigma d^{\alpha /2}$, $u=\log
(z_{2}/z_{1})/a$ and $v\left( u\right) =a/2+u$. It is well known that the
joint probability distribution function of $\left( \eta \left( x_{1}\right)
,\eta \left( x_{2}\right) \right) $ is given by (see e.g. Huser and Davison
(2013)) 
\begin{equation*}
F_{x_{1},x_{2}}\left( z_{1},z_{2}\right) =\mathbb{P}\left[ \eta \left(
x_{1}\right) \leq z_{1},\eta \left( x_{2}\right) \leq z_{2}\right] =\exp
\left( -V_{x_{1},x_{2}}\left( z_{1},z_{2}\right) \right) ,
\end{equation*}%
where%
\begin{equation*}
V_{x_{1},x_{2}}\left( z_{1},z_{2}\right) =\frac{1}{z_{1}}\Phi \left( v\left(
u\right) \right) +\frac{1}{z_{2}}\Phi \left( v\left( -u\right) \right)
,\qquad z_{1},z_{2}>0.  \label{Eq_V_dim_2}
\end{equation*}%
Here $\Phi $ denotes the cumulative distribution function of the standard
Gaussian distribution. The term $V_{x_{1},x_{2}}$ is referred to as the
pairwise exponent function. Let us now consider the \textquotedblleft
normalized\textquotedblright\ (linear) increments of the logarithm of the
Brown-Resnick random field%
\begin{equation*}
U=d^{-\alpha /2}\sigma ^{-1}\log \left( \eta \left( x_{2}\right) /\eta
\left( x_{1}\right) \right) .
\end{equation*}%
The following proposition provides the conditional and marginal
distributions of $U$ and allows us to deduce that it has asymptotically a
standard Gaussian distribution as the distance $d$ tends to $0$. Such a
result generalizes Proposition 3 in Robert (2020).

\begin{proposition}
\label{Prop_cond_dist_U}The conditional distribution of $U$ given $\eta
\left( x_{1}\right) =\eta >0$ is characterized by%
\begin{equation*}
\mathbb{P}\left[ \left. U\leq u\right\vert \eta \left( x_{1}\right) =\eta %
\right] =\exp \left( -\frac{1}{\eta }\left[ V_{x_{1},x_{2}}(1,e^{\sigma
d^{\alpha /2}u})-1\right] \right) \Phi \left( v\left( u\right) \right)
,\qquad u\in \mathbf{R,}
\end{equation*}%
and its marginal distribution by%
\begin{equation*}
\mathbb{P}\left[ U\leq u\right] =\frac{\Phi \left( v\left( u\right) \right) 
}{V_{x_{1},x_{2}}(1,e^{\sigma d^{\alpha /2}u})},\qquad u\in \mathbf{R}.
\end{equation*}%
It follows that 
\begin{equation*}
\lim_{d\rightarrow 0}\mathbb{P}\left[ U\leq u\right] =\Phi \left( u\right)
,\qquad u\in \mathbf{R}.
\end{equation*}
\end{proposition}

The fact that the asymptotic distribution of $U$ (as $d$ tends to $0$) is a
standard Gaussian distribution is not a surprise since the probability that $%
x_{1}$ and $x_{2}$ belong to the same cell of the canonical tessellation of
the max-stable random field tends to $1$. Indeed, in a common cell, the
values of the max-stable random field are generated by the same isotropic
fractional Brownian random field. It is natural to first
study the asymptotic behaviors of the increment sums for an isotropic
fractional Brownian random field before considering a Brown-Resnick random
field (see Section \ref{Sec_Asymp_fBf}).

The distribution of $\left( \eta \left( x_{1}\right) ,\eta \left(
x_{2}\right) \right) $ is absolutely continuous with respect to the Lebesgue
measure on $(\mathbf{R}_{+})^{2}$. Its density function satisfies%
\begin{equation*}
f_{x_{1},x_{2}}\left( z_{1},z_{2}\right) =\frac{\partial }{\partial
z_{1}\partial z_{2}}F_{x_{1},x_{2}}\left( z_{1},z_{2}\right) ,\qquad
z_{1},z_{2}\in \mathbf{R}_{+},
\end{equation*}%
and will be used for the contribution of the pair $\left( \eta \left(
x_{1}\right) ,\eta \left( x_{2}\right) \right) $ to the pairwise CL
function. For any $\alpha \in (0,2)$, $\sigma ^{2}>0$ and $z_{1},z_{2}\in 
\mathbf{R}_{+}$, this joint density function is a differentiable function
with respect to $\left( \alpha ,\sigma \right) $. The following proposition
provides the asymptotic contributions of the pair to the pairwise score
functions.

\begin{proposition}
\label{Prop_pairwise_pdf}Let $u\in \mathbf{R}$ be fixed. Let $x_{1},x_{2}\in 
\mathbf{R}^{2}$ and $z_{1},z_{2}\in \mathbf{R}_{+}$ be such that $d^{-\alpha
/2}\sigma ^{-1}\log \left( z_{2}/z_{1}\right) =u$, where $d=\left\Vert
x_{2}-x_{1}\right\Vert >0$. Then 
\begin{eqnarray*}
\lim_{d\rightarrow 0}\frac{\partial }{\partial \sigma }\log
f_{x_{1},x_{2}}\left( z_{1},z_{2}\right) &=&\frac{1}{\sigma }\left(
u^{2}-1\right) , \\
\lim_{d\rightarrow 0}\frac{1}{\log d}\frac{\partial }{\partial \alpha }\log
f_{x_{1},x_{2}}\left( z_{1},z_{2}\right) &=&\frac{1}{2}\left( u^{2}-1\right)
.
\end{eqnarray*}
\end{proposition}

The asymptotic score contributions of a pair are therefore proportional to $%
\left( u^{2}-1\right) $. Further, $u$ will be replaced by the normalized
increment of $\log \left( \eta \right) $ which has asymptotically a standard
Gaussian distribution as stated in Proposition \ref{Prop_cond_dist_U}. This
fact ensures that the asymptotic score contributions are asymptotically
unbiased.

\subsection{Triplewise joint distributions and asymptotic score contributions%
}

Let us now consider three sites $x_{1},x_{2},x_{3}\in \mathbf{R}^{2}$ and
denote by $d_{1,2}=\left\Vert x_{2}-x_{1}\right\Vert $, $d_{1,3}=\left\Vert
x_{3}-x_{1}\right\Vert $, $d_{2,3}=\left\Vert x_{3}-x_{2}\right\Vert $ the
distances between two different sites. Let $z_{1},z_{2},z_{3}\in \mathbf{R}%
_{+}$ and, for $i,j=1,2,3$ such that $i\neq j$, let $a_{i,j}=\sigma
d_{i,j}^{\alpha /2}$, $u_{i,j}=\log (z_{j}/z_{i})/a_{i,j}$ and $%
v_{i,j}\left( u\right) =a_{i,j}/2+u_{i,j}$. The joint probability
distribution function of $\left( \eta \left( x_{1}\right) ,\eta \left(
x_{2}\right) ,\eta \left( x_{3}\right) \right) $ is given by (see e.g. Huser
and Davison (2013)) 
\begin{equation*}
F_{x_{1},x_{2},x_{3}}\left( z_{1},z_{2},z_{3}\right) =\mathbb{P}\left[ \eta
\left( x_{1}\right) \leq z_{1},\eta \left( x_{2}\right) \leq z_{2},\eta
\left( x_{3}\right) \leq z_{3}\right] =\exp \left(
-V_{x_{1},x_{2},x_{3}}\left( z_{1},z_{2},z_{3}\right) \right) ,
\end{equation*}%
where%
\begin{multline}
 V_{x_{1},x_{2},x_{3}}\left( z_{1},z_{2},z_{3}\right)   = \frac{1}{z_{1}}\Phi _{2}\left( \left( 
\begin{array}{c}
v_{1,2}\left( u_{1,2}\right) \\ 
v_{1,3}\left( u_{1,3}\right)%
\end{array}%
\right) ;\left( 
\begin{array}{cc}
1 & R_{1} \\ 
R_{1} & 1%
\end{array}%
\right) \right) + \frac{1}{z_{2}}\Phi _{2}\left( \left( 
\begin{array}{c}
v_{1,2}\left( -u_{1,2}\right) \\ 
v_{2,3}\left( u_{2,3}\right)%
\end{array}%
\right) ;\left( 
\begin{array}{cc}
1 & R_{2} \\ 
R_{2} & 1%
\end{array}%
\right) \right)\\
 +\frac{1}{z_{3}}\Phi _{2}\left( \left( 
\begin{array}{c}
v_{1,3}\left( -u_{1,3}\right) \\ 
v_{2,3}\left( -u_{2,3}\right)%
\end{array}%
\right) ;\left( 
\begin{array}{cc}
1 & R_{3} \\ 
R_{3} & 1%
\end{array}%
\right) \right)  \label{Eq_V_triple}
\end{multline}
with%
\begin{equation*}
R_{1}=\frac{d_{1,2}^{\alpha }+d_{1,3}^{\alpha }-d_{2,3}^{\alpha }}{2\left(
d_{1,2}d_{1,3}\right) ^{\alpha /2}},\qquad R_{2}=\frac{d_{1,2}^{\alpha
}+d_{2,3}^{\alpha }-d_{1,3}^{\alpha }}{2\left( d_{1,2}d_{2,3}\right)
^{\alpha /2}},\qquad R_{3}=\frac{d_{1,3}^{\alpha }+d_{2,3}^{\alpha
}-d_{1,2}^{\alpha }}{2\left( d_{1,3}d_{2,3}\right) ^{\alpha /2}}.
\label{q_Corr}
\end{equation*}%
Here $\Phi _{2}\left( \bullet ,\Sigma \right) $ denotes the bivariate
cumulative distribution function of the centered Gaussian distribution with
covariance matrix $\Sigma $. As for the pairs, let us also consider the
\textquotedblleft normalized\textquotedblright\ (linear) increments of the
logarithm of the Brown-Resnick random field%
\begin{equation*}
U_{1,2}=d_{1,2}^{-\alpha /2}\sigma ^{-1}\log \left( \eta \left( x_{2}\right)
/\eta \left( x_{1}\right) \right) ,\qquad U_{1,3}=d_{1,3}^{-\alpha /2}\sigma
^{-1}\log \left( \eta \left( x_{3}\right) /\eta \left( x_{1}\right) \right) .
\end{equation*}%
The following proposition provides the conditional and marginal
distributions of the vector $\left( U_{1,2},U_{1,3}\right) $ and allows us
to deduce that it has asymptotically a bivariate Gaussian distribution as
the distances $d_{1,2}$ and $d_{1,3}$ tend to $0$ proportionally.

\begin{proposition}
\label{Prop_cond_dist_U_1_U_2}The conditional distribution of $%
(U_{1,2},U_{1,3})$ given $\eta \left( x_{1}\right) =\eta >0$ is
characterized by%
\begin{multline*}
\mathbb{P}\left[ {\left. U_{1,2}\leq u_{2},U_{1,3}\leq u_{3}\right\vert }%
\eta \left( x_{1}\right) =\eta \right]  = \exp \left( -\frac{1}{\eta }\left[ V_{x_{1},x_{2},x_{3}}(1,e^{\sigma
d_{12}^{\alpha /2}u_{2}},e^{\sigma d_{13}^{\alpha /2}u_{3}})-1\right]
\right)\\
\times  \Phi _{2}\left( 
\begin{pmatrix}
v_{1,2}\left( u_{2}\right) \\ 
v_{1,3}\left( u_{3}\right)%
\end{pmatrix}%
;%
\begin{pmatrix}
1 & R_{1} \\ 
R_{1} & 1%
\end{pmatrix}%
\right) ,\qquad u_{1},u_{2}\in \mathbf{R},
\end{multline*}%
and its marginal distribution by 
\begin{equation*}
\mathbb{P}\left[ U_{1,2}\leq u_{2},U_{1,3}\leq u_{3}\right] =\frac{\Phi
_{2}\left( 
\begin{pmatrix}
v_{1,2}\left( u_{2}\right) \\ 
v_{1,3}\left( u_{3}\right)%
\end{pmatrix}%
;%
\begin{pmatrix}
1 & R_{1} \\ 
R_{1} & 1%
\end{pmatrix}%
\right) }{V_{x_{1},x_{2},x_{3}}(1,e^{\sigma d_{12}^{\alpha
/2}u_{2}},e^{\sigma d_{13}^{\alpha /2}u_{3}})},\qquad u_{1},u_{2}\in \mathbf{%
R}.
\end{equation*}%
It follows that, if $\left\Vert x_{2}-x_{1}\right\Vert =\delta d_{1,2}$, $%
\left\Vert x_{3}-x_{1}\right\Vert =\delta d_{1,3}$, $\left\Vert
x_{3}-x_{2}\right\Vert =\delta d_{2,3}$, where $d_{i,j}$, $i\neq j$, is
fixed, then%
\begin{equation*}
\lim_{\delta \rightarrow 0}\mathbb{P}\left[ U_{1,2}\leq u_{2},U_{1,3}\leq
u_{3}\right] =\Phi _{2}\left( 
\begin{pmatrix}
u_{2} \\ 
u_{3}%
\end{pmatrix}%
;%
\begin{pmatrix}
1 & R_{1} \\ 
R_{1} & 1%
\end{pmatrix}%
\right) ,\qquad u_{1},u_{2}\in \mathbf{R}.
\end{equation*}
\end{proposition}

The comment concerning the asymptotic distribution of $U$ also holds for $%
\left( U_{1,2},U_{1,3}\right) $. The probability that $x_{1}$, $x_{2}$ and $%
x_{3}$ belong to the same cell of the canonical tessellation of the
max-stable random field tends to $1$ as $\delta $ tends to $0$. Therefore
the vector $\left( U_{1,2},U_{1,3}\right) $ tends to have the same
distribution as the vector of normalized linear increments of an isotropic
fractional Brownian random field.

The distribution of $\left( \eta \left( x_{1}\right) ,\eta \left(
x_{2}\right) ,\eta \left( x_{3}\right) \right) $ is absolutely continuous
with respect to the Lebesgue measure on $(\mathbf{R}_{+})^{3}$. Its density
function satisfies%
\begin{equation*}
f_{x_{1},x_{2},x_{3}}\left( z_{1},z_{2},z_{3}\right) =\frac{\partial }{%
\partial z_{1}\partial z_{2}\partial z_{3}}F_{x_{1},x_{2},x_{3}}\left(
z_{1},z_{2},z_{3}\right) ,\qquad z_{1},z_{2},z_{3}\in \mathbf{R}_{+},
\end{equation*}%
and will be used for the contribution of the triple $\left( \eta \left(
x_{1}\right) ,\eta \left( x_{2}\right) ,\eta \left( x_{3}\right) \right) $
to the triplewise CL function. For any $\alpha \in (0,2)$, $\sigma ^{2}>0$
and $z_{1},z_{2},z_{3}\in \mathbf{R}_{+}$, this joint density function is a
differentiable function with respect to $\left( \alpha ,\sigma \right) $.
The following proposition provides the asymptotic contributions of the
triple $\left( \eta \left( x_{1}\right) ,\eta \left( x_{2}\right) ,\eta
\left( x_{3}\right) \right) $ to the triplewise score functions.

\begin{proposition}
\label{Prop_triplewise_pdf}Let $u_{2},u_{3}\in \mathbf{R}$ be fixed. Let $%
x_{1},x_{2},x_{3}\in \mathbf{R}^{2}$ and $z_{1},z_{2},z_{3}\in \mathbf{R}%
_{+} $ be such that%
\begin{equation*}
\delta ^{-\alpha /2}d_{1,2}^{-\alpha /2}\sigma ^{-1}\log \left(
z_{2}/z_{1}\right) =u_{2}\text{ and }\delta ^{-\alpha /2}d_{1,3}^{-\alpha
/2}\sigma ^{-1}\log \left( z_{3}/z_{1}\right) =u_{3},
\end{equation*}%
where $\delta d_{1,2}=\left\Vert x_{2}-x_{1}\right\Vert $, $\delta
d_{1,3}=\left\Vert x_{3}-x_{1}\right\Vert $, $\delta d_{2,3}=\left\Vert
x_{3}-x_{2}\right\Vert $. Then 
\begin{eqnarray*}
\lim_{\delta \rightarrow 0}\frac{\partial }{\partial \sigma }\log
f_{x_{1},x_{2},x_{3}}\left( z_{1},z_{2},z_{3}\right) &=&\frac{1}{\sigma }%
\left( \left( 
\begin{array}{cc}
u_{2} & u_{3}%
\end{array}%
\right) \left( 
\begin{array}{cc}
1 & R_{1} \\ 
R_{1} & 1%
\end{array}%
\right) ^{-1}\left( 
\begin{array}{c}
u_{2} \\ 
u_{3}%
\end{array}%
\right) -2\right) , \\
\lim_{\delta \rightarrow 0}\frac{1}{\log \delta }\frac{\partial }{\partial
\alpha }\log f_{x_{1},x_{2},x_{3}}\left( z_{1},z_{2},z_{3}\right) &=&\frac{1%
}{2}\left( \left( 
\begin{array}{cc}
u_{2} & u_{3}%
\end{array}%
\right) \left( 
\begin{array}{cc}
1 & R_{1} \\ 
R_{1} & 1%
\end{array}%
\right) ^{-1}\left( 
\begin{array}{c}
u_{2} \\ 
u_{3}%
\end{array}%
\right) -2\right) .
\end{eqnarray*}
\end{proposition}

The asymptotic contributions of a triple are therefore proportional to a
quadratic function of $\left( u_{2},u_{3}\right) $. Further, $u_{2}$ and $%
u_{3}$ will be replaced by normalized increments of $\log \left( \eta
\right) $ over a triangle of the Delaunay triangulation (see Section \ref%
{Sec_randomized_sampling_scheme}) which have asymptotically a bivariate
standard Gaussian distribution with correlation coefficient $R_{1}$, as
stated in Proposition \ref{Prop_cond_dist_U_1_U_2}. We can also conclude
that the asymptotic score contributions are asymptotically unbiased.

If we let $\delta ^{-\alpha /2}d_{1,2}^{-\alpha /2}\sigma ^{-1}\log \left(
z_{1}/z_{2}\right) =\tilde{u}_{1}$ and $\delta ^{-\alpha /2}d_{2,3}^{-\alpha
/2}\sigma ^{-1}\log \left( z_{3}/z_{2}\right) =\tilde{u}_{3}$ with fixed $%
\tilde{u}_{1},\tilde{u}_{3}\in \mathbf{R}$, we also get%
\begin{equation*}
\lim_{\delta \rightarrow 0}\frac{\partial }{\partial \sigma }\log
f_{x_{1},x_{2},x_{3}}\left( z_{1},z_{2},z_{3}\right) =\frac{1}{\sigma }%
\left( \left( 
\begin{array}{cc}
\tilde{u}_{1} & \tilde{u}_{3}%
\end{array}%
\right) \left( 
\begin{array}{cc}
1 & R_{2} \\ 
R_{2} & 1%
\end{array}%
\right) ^{-1}\left( 
\begin{array}{c}
\tilde{u}_{1} \\ 
\tilde{u}_{3}%
\end{array}%
\right) -2\right) .
\end{equation*}%
In particular, an invariance property with respect to the choice of the
order of the points $x_{1}$, $x_{2}$ and $x_{3}$ holds for the triplewise
score functions. However, it will be necessary to order these points later.

\section{The weighted CL approach}

\subsection{The randomized sampling scheme\label%
{Sec_randomized_sampling_scheme}}

We assume that the data-sites are given by a realization of a homogeneous
Poisson point process of intensity $N$ in $\mathbf{R}^{2}$, denoted by $%
P_{N} $, which is independent of the Brown-Resnick random field. Let us
denote by $\mathbf{C}=(-1/2,1/2]^{2}$ the square where we will consider the
sites for the observations of the max-stable random field.

The Delaunay graph Del$(P_{N})$ based on $P_{N}$ is our connection scheme
and is defined as the unique triangulation with vertices in $P_{N}$ such
that the circumball of each triangle contains no point of $P_{N}$ in its
interior. With a slight abuse of notation, we identify $\text{Del}(P_{N})$
to its skeleton. When $x_{1},x_{2}\in P_{N}$ are Delaunay neighbors, we
write $x_{1}\sim x_{2}$ in $\text{Del}(P_{N})$.

For a Borel subset $\mathbf{B}$ in $\mathbf{R}^{2}$, let $E_{N,\mathbf{B}}$
be\ the set of couples $(x_{1},x_{2})$ such that the following conditions
hold:%
\begin{equation*}
x_{1}\sim x_{2}\text{ in Del}(P_{N}),\quad x_{1}\in \mathbf{B},\quad \text{%
and}\quad x_{1}\preceq x_{2},
\end{equation*}%
where $\preceq $ denotes the lexicographic order. When $\mathbf{B}=\mathbf{C}
$, we only write $E_{N}=E_{N,\mathbf{C}}$.

For a Borel subset $\mathbf{B}$ in $\mathbf{R}^{2}$, let $DT_{N,\mathbf{B}}$
be the set of triples $(x_{1},x_{2},x_{3})$ satisfying the following
properties%
\begin{equation*}
\Delta (x_{1},x_{2},x_{3})\in \text{Del}(P_{N}),\quad x_{1}\in \mathbf{B}%
,\quad \text{and}\quad x_{1}\preceq x_{2}\preceq x_{3},
\end{equation*}%
where $\Delta (x_{1},x_{2},x_{3})$ is the convex hull of $%
(x_{1},x_{2},x_{3}) $. When $\mathbf{B}=\mathbf{C}$, we only write $%
DT_{N}=DT_{N,\mathbf{C}}$.

\subsection{The weighted CL objective functions and the CL estimators}

The (tapered) pairwise CL objective function\ is defined as 
\begin{equation*}
\ell _{2,N}\left( \sigma ,\alpha \right) =\sum_{\left( x_{1},x_{2}\right)
\in E_{N}}\log f_{x_{1},x_{2}}\left( \eta \left( x_{1}\right) ,\eta \left(
x_{2}\right) \right) ,
\end{equation*}%
while the (tapered) triplewise CL objective function\ is defined as%
\begin{equation*}
\ell _{3,N}\left( \sigma ,\alpha \right) =\sum_{\left(
x_{1},x_{2},x_{3}\right) \in DT_{N}}\log f_{x_{1},x_{2},x_{3}}\left( \eta
\left( x_{1}\right) ,\eta \left( x_{2}\right) ,\eta \left( x_{3}\right)
\right) .
\end{equation*}%
Thereby, in the CL objective functions, we exclude pairs that are not edges
of the Delaunay triangulation or triples that are not vertices of triangles
of this triangulation. Restricting the CL objective functions to the most
informative pairs and triples for the estimation of the parameters does not
modify the approach that follows, but allows us to simplify the presentation
and the proofs.

From Section 4.4 of Dombry et al. (2018), we know that there exist families
of positive functions $\left( l_{x_{1},x_{2}}\right) _{x_{1},x_{2}\in 
\mathbf{R}^{2}}$ and $\left( l_{x_{1},x_{2},x_{3}}\right)
_{x_{1},x_{2},x_{3}\in \mathbf{R}^{2}}$ with $l_{x_{1},x_{2}}:\mathbf{R}%
^{2}\rightarrow \mathbf{R}$ and $l_{x_{1},x_{2},x_{3}}:\mathbf{R}%
^{3}\rightarrow \mathbf{R}$ such that the following Lipschitz conditions
hold: for any $\sigma _{1},\sigma _{2}>0$ and $\alpha _{1},\alpha _{2}\in
(0,2)$ 
\begin{equation*}
\left\vert \log \frac{f_{x_{1},x_{2}}\left( z_{1},z_{2};\left( \sigma
_{2},\alpha _{2}\right) \right) }{f_{x_{1},x_{2}}\left( z_{1},z_{2};\left(
\sigma _{1},\alpha _{1}\right) \right) }\right\vert \leq
l_{x_{1},x_{2}}\left( z_{1},z_{2}\right) \left( \left\vert \sigma
_{2}-\sigma _{1}\right\vert +\left\vert \alpha _{2}-\alpha _{1}\right\vert
\right)
\end{equation*}%
and%
\begin{equation*}
\left\vert \log \frac{f_{x_{1},x_{2},x_{3}}\left( z_{1},z_{2},z_{3};\left(
\sigma _{2},\alpha _{2}\right) \right) }{f_{x_{1},x_{2},x_{3}}\left(
z_{1},z_{2},z_{3};\left( \sigma _{1},\alpha _{1}\right) \right) }\right\vert
\leq l_{x_{1},x_{2},x_{3}}\left( z_{1},z_{2,z_{3}}\right) \left( \left\vert
\sigma _{2}-\sigma _{1}\right\vert +\left\vert \alpha _{2}-\alpha
_{1}\right\vert \right) .
\end{equation*}

Let us denote by $\left( \sigma _{0},\alpha _{0}\right) $ the true
parameters. We assume that $\sigma _{0}$ belongs to a compact set $\mathrm{S}%
_{\sigma }$ of $\mathbf{R}_{+}$ and that $\alpha _{0}$ belongs to a compact
set $\mathrm{S}_{\alpha }$ of $(0,2)$. We can now define the MCLEs of $%
\sigma $ and $\alpha $.

When $\alpha _{0}$ is assumed to be known, the pairwise and triplewise
maximum (tapered) CL estimators of $\sigma _{0}$, $\hat{\sigma}_{j,N}$, are
respectively defined as a solution of the maximization problems%
\begin{equation*}
\max_{\sigma \in \mathrm{S}_{\sigma }}\ell _{j,N}\left( \sigma ,\alpha
_{0}\right) ,\qquad j=2,3.
\end{equation*}%
When $\sigma _{0}$ is assumed to be known, the pairwise and triplewise
maximum (tapered) CL estimators of $\alpha _{0}$, $\hat{\alpha}_{j,N}$, are
respectively defined as a solution of the maximization problems%
\begin{equation*}
\max_{\alpha \in \mathrm{S}_{\alpha }}\ell _{j,N}\left( \sigma _{0},\alpha
\right) ,\qquad j=2,3.
\end{equation*}%
Note that the solutions of these maximization problems become unique as $%
N\rightarrow \infty $. This can be viewed from the first-order optimality
conditions and the asymptotic approximations of the score functions obtained
in Propositions \ref{Prop_pairwise_pdf} and \ref{Prop_triplewise_pdf}.

\section{Main results}

Our aim is to characterize the asymptotic distributions of the MCLEs. We
provide intermediate results for different random fields in order to
understand how we obtained the different families of asymptotic
distributions of our estimators. We first provide some definitions and
notations related to the Poisson-Delaunay triangulation. Then we consider
sums of square increments of an isotropic fractional Brownian field on the
edges of the Delaunay triangles and provide Central Limit Theorems using
Malliavin calculus. We only consider the case $\alpha \in (0,1)$ for which
the asymptotic distributions are Gaussian. This is not a very restrictive
constraint since almost all empirical studies that use the spatial
Brown-Resnick random field obtain values for $\alpha $ in this interval (see
e.g. Davison et al. (2012), Engelke et al. (2014), Einmahl et al. (2015) or
de Fondeville and Davison (2018)). Third we consider sums of square
increments of the pointwise maximum of two independent isotropic fractional
Brownian fields and show that the asymptotic behaviors of the sums now
depend on the local time at the level $0$ of the difference between the two
fractional Brownian fields. Fourth we generalize these results to the
max-stable Brown-Resnick random field and, using approximation of the
pairwise and triplewise CL objective functions, we derive the asymptotic
properties of the MCLEs.

\subsection{Definitions and notations}

A classical object in Stochastic Geometry is the typical cell. To define it,
let us consider a Delaunay triangulation Del$(P_{1})$ based on a homogeneous
Poisson point process of intensity $1$. With each cell $C\in \text{Del}%
(P_{1})$, we associate the circumcenter $z(C)$ of $C$. Now, let $\mathbf{B}$
be a Borel subset in $\mathbf{R}^{2}$ with area $a(\mathbf{B})\in (0,\infty
) $. The cell intensity $\beta _{2}$ of $\text{Del}(P_{1})$ is defined as
the mean number of cells per unit area, i.e. 
\begin{equation*}
\beta _{2}=\frac{1}{a(\mathbf{B})}\mathbb{E}\left[ {|\{C\in \text{Del}%
(P_{1}):z(C)\in \mathbf{B}\}|}\right] .
\end{equation*}%
It is known that $\beta _{2}=2$, see e.g. Theorem 10.2.9. in Schneider and
Weil (2008). Then, we define the typical cell as a random triangle $\mathcal{%
C}$ with distribution given as follows: for any positive measurable and
translation invariant function $g:\mathcal{K}_{2}\rightarrow \mathbf{R}$, we
have 
\begin{equation*}
\mathbb{E}\left[ {g(\mathcal{C})}\right] =\frac{1}{\beta _{2}a(\mathbf{B})}%
\mathbb{E}\left[ {\sum_{C\in \text{Del}(P_{1}):z(C)\in \mathbf{B}}g(C)}%
\right] ,
\end{equation*}%
where $\mathcal{K}_{2}$ denotes the set of convex compact subsets in $%
\mathbf{R}^{2}$, endowed with the Fell topology (see Section 12.2 in
Schneider and Weil (2008) for the definition). The distribution of $\mathcal{%
C}$ has the following integral representation (see e.g. Theorem 10.4.4. in
Schneider and Weil (2008)): 
\begin{equation}
\mathbb{E}\left[ {g(\mathcal{C})}\right] =\frac{1}{6}\int_{0}^{\infty
}\int_{(\mathbf{S}^{1})^{3}}r^{3}e^{-\pi r^{2}}a(\Delta
(u_{1},u_{2},u_{3}))g(\Delta (ru_{1},ru_{2},ru_{3}))\sigma (\mathrm{d}%
u_{1})\sigma (\mathrm{d}u_{2})\sigma (\mathrm{d}u_{3})\mathrm{d}r,
\label{eq:typicalcell}
\end{equation}%
where $\mathbf{S}^{1}$ is the unit sphere of $\mathbf{R}^{2}$ and $\sigma $
is the spherical Lebesgue measure on $\mathbf{S}^{1}$ with normalization $%
\sigma \left( \mathbf{S}^{1}\right) =2\pi $. It means that $\mathcal{C}$ is
equal in distribution to $R\Delta (U_{1},U_{2},U_{3})$, where $R$ and $%
(U_{1},U_{2},U_{3})$ are independent with probability density functions
given respectively by $2\pi ^{2}r^{3}e^{-\pi r^{2}}$ and $a(\Delta
(u_{1},u_{2},u_{3}))/(12\pi ^{2})$.

In a similar way, we can define the notion of typical edge. The edge
intensity $\beta _{1}$ of $\text{Del}(P_{1})$ is defined as the mean number
of edges per unit area and is equal to $\beta _{1}=3$ (see e.g. Theorem
10.2.9. in Schneider and Weil (2008)). The distribution of the length of the
typical edge is the same as the distribution of $D=R||U_{1}-U_{2}||$. Its
probability density function $f_{D}$ satisfies the following equality 
\begin{equation}
\mathbb{P}\left[ D\leq \ell\right] =\int_{0}^{\ell}f_{D}(d)\mathrm{d}d=\frac{\pi }{%
3}\int_{0}^{\infty }\int_{(\mathbf{S}^{1})^{2}}r^{3}e^{-\pi r^{2}}a(\Delta
(u_{1},u_{2},e_{1}))\mathbb{I}\left[ r\left\Vert u_{1}-u_{2}\right\Vert \leq
\ell\right] \sigma (\mathrm{d}u_{1})\sigma (\mathrm{d}u_{2})\mathrm{d}r,
\label{eq:typicallength}
\end{equation}%
where $e_{1}=(1,0)$ and $\ell>0$. Following Eq. $\left( \ref{eq:typicalcell}%
\right) $, a typical couple of (distinct) Delaunay edges with a common
vertex can be defined as a $3$-tuple of random variables $%
(D_{1},D_{2},\Theta )$, where $D_{1},D_{2}\geq 0$ and $\Theta \in \lbrack -%
\frac{\pi }{2},\frac{\pi }{2})$, with distribution given by 
\begin{multline*}
\mathbb{P}[(D_{1},D_{2},\Theta )\in B]=\frac{1}{6}\int_{0}^{\infty }\int_{(%
\mathbf{S}^{1})^{3}}r^{3}e^{-\pi r^{2}}a(\Delta (u_{1},u_{2},u_{3})) \\
\times \mathbb{I}[(r||u_{3}-u_{2}||,r||u_{2}-u_{1}||,\arcsin \left( \cos
(\zeta _{u_{1},u_{2}}/2)\right) )\in B]\sigma (\mathrm{d}u_{1})\sigma (%
\mathrm{d}u_{2})\sigma (\mathrm{d}u_{3})\mathrm{d}r,
\end{multline*}%
where $\zeta _{u_{1},u_{2}}$ is the measure of the angle $(u_{1},u_{2})$ and
where $B$ is any Borel subset in $\mathbf{R}_{+}^{2}\times \lbrack -\frac{%
\pi }{2},\frac{\pi }{2})$. The random variables $D_{1},D_{2}$ (resp. $\Theta 
$) can be interpreted as the lengths of the two typical edges (resp. as the
angle between the edges). In particular, the length of a typical edge is
equal in distribution to $D=R||U_{2}-U_{1}||$ with distribution given in Eq. 
$\left( \ref{eq:typicallength}\right) $.

\subsection{Asymptotic distributions of squared increment sums for an
isotropic fractional Brownian field\label{Sec_Asymp_fBf}}

Let $\left( W\left( x\right) \right) _{x\in \mathbf{R}^{2}}$ be an isotropic
fractional Brownian field where $W\left( 0\right) =0$ a.s. and $\text{var}\left(
W\left( x\right) \right) =\sigma ^{2}\left\Vert x\right\Vert ^{\alpha }$ for
some $\alpha \in (0,1)$ and $\sigma ^{2}>0$. For two sites $x_{1},x_{2}\in 
\mathbf{R}^{2}$, let us define the normalized increment between $x_{1}$ and $%
x_{2}$ as%
\begin{equation*}
U_{x_{1},x_{2}}^{(W)}=\sigma ^{-1}d_{1,2}^{-\alpha /2}\left( W\left(
x_{2}\right) -W\left( x_{1}\right) \right)
\end{equation*}%
with $d_{1,2}=\left\Vert x_{2}-x_{1}\right\Vert $.

The (normalized) squared increment sum for the edges of the Delaunay
triangulation is given by%
\begin{equation*}
V_{2,N}^{(W)}=\frac{1}{\sqrt{\left\vert E_{N}\right\vert }}\sum_{\left(
x_{1},x_{2}\right) \in E_{N}}\left( (U_{x_{1},x_{2}}^{(W)})^{2}-1\right) ,
\end{equation*}%
while the (normalized) squared increment sum for the pairs of edges of
Delaunay triangles is defined as%
\begin{equation*}
V_{3,N}^{(W)}=\frac{1}{\sqrt{\left\vert DT_{N}\right\vert }}\sum_{\left(
x_{1},x_{2},x_{3}\right) \in DT_{N}}\left( \left( 
\begin{array}{cc}
U_{x_{1},x_{2}}^{(W)} & U_{x_{1},x_{3}}^{(W)}%
\end{array}%
\right) \left( 
\begin{array}{cc}
1 & R_{x_{1},x_{2},x_{3}} \\ 
R_{x_{1},x_{2},x_{3}} & 1%
\end{array}%
\right) ^{-1}\left( 
\begin{array}{c}
U_{x_{1},x_{2}}^{(W)} \\ 
U_{x_{1},x_{3}}^{(W)}%
\end{array}%
\right) -2\right) ,
\end{equation*}%
where%
\begin{equation}
R_{x_{1},x_{2},x_{3}}=\text{corr}(U_{x_{1},x_{2}}^{(W)},U_{x_{1},x_{3}}^{(W)})=%
\frac{d_{1,2}^{\alpha }+d_{1,3}^{\alpha }-d_{2,3}^{\alpha }}{2\left(
d_{1,2}d_{1,3}\right) ^{\alpha /2}},  \label{eq:coorU}
\end{equation}%
with $d_{1,3}=\left\Vert x_{3}-x_{1}\right\Vert >0$ and $d_{2,3}=\left\Vert
x_{3}-x_{2}\right\Vert >0$. Let%
\begin{equation*}
\tilde{U}_{x_{1},x_{2},x_{3}}^{(W)}=(1-R_{x_{1},x_{2},x_{3}}^{2})^{-1/2}%
\left(
U_{x_{1},x_{2}}^{(W)}-R_{x_{1},x_{2},x_{3}}U_{x_{1},x_{3}}^{(W)}\right)
\quad \text{and}\quad \tilde{U}_{x_{1},x_{3}}^{(W)}=U_{x_{1},x_{3}}^{(W)}.
\end{equation*}%
Note that $\tilde{U}_{x_{1},x_{2},x_{3}}^{(W)}$ is a normalized increment
based on the three points $x_{1},x_{2},x_{3}$ (see e.g. Chan and Wood
(2002)) and that%
\begin{equation*}
\text{corr}(\tilde{U}_{x_{1},x_{2},x_{3}}^{(W)},\tilde{U}_{x_{1},x_{3}}^{(W)})=0%
\text{.}
\end{equation*}%
The sum $V_{3,N}^{(W)}$ may be rewritten as%
\begin{equation*}
V_{3,N}^{(W)}=\frac{1}{\sqrt{\left\vert DT_{N}\right\vert }}\sum_{\left(
x_{1},x_{2},x_{3}\right) \in DT_{N}}\left( [(\tilde{U}%
_{x_{1},x_{2},x_{3}}^{(W)})^{2}-1]+[(\tilde{U}_{x_{1},x_{3}}^{(W)})^{2}-1]%
\right) .
\end{equation*}

The following theorem states that the asymptotic distributions of $%
V_{2,N}^{(W)}$ and $V_{3,N}^{(W)}$ are Gaussian. Their asymptotic variances
are known, but quite intricate. We provide their integral representations in
Section 1  in the Supplementary Material. %\ref{Sect_Int_Rep_variance}

\begin{theorem}
\label{Prop_Conv_V_n_Gaussian}Let $\alpha \in (0,1)$. Then there exist constants $%
\sigma _{V_{2}}^{2}>0$ and $\sigma _{V_{3}}^{2}>0$ such that, as $%
N\rightarrow \infty $,%
\begin{equation*}
V_{2,N}^{(W)}\overset{\mathcal{D}}{\rightarrow }\mathcal{N}\left( 0,\sigma
_{V_{2}}^{2}\right) ,\qquad V_{3,N}^{(W)}\overset{\mathcal{D}}{\rightarrow }%
\mathcal{N}\left( 0,\sigma _{V_{3}}^{2}\right) .
\end{equation*}
\end{theorem}

We note that the rates of convergences of both sums are the same as in
Theorem 3.2\ of Chan and Wood (2000) or in Theorem 1 of Zhu and Stein (2002)
where statistics based on square increments on regular grids have been
considered.

\subsection{Asymptotic distributions of squared increment sums for the
(pointwise) maximum of two independent fractional Brownian fields}

Let $\left( W^{(1)}\left( x\right) \right) _{x\in \mathbf{R}^{2}}$ and $%
\left( W^{(2)}\left( x\right) \right) _{x\in \mathbf{R}^{2}}$ be two
independent isotropic fractional Brownian fields, where $W^{(1)}\left(
0\right) =W^{(2)}\left( 0\right) =0$ a.s. and $\text{var}\left( W^{(1)}\left(
x\right) \right) =\text{var}\left( W^{(2)}\left( x\right) \right) =\sigma
^{2}\left\Vert x\right\Vert ^{\alpha }$ for some $\alpha \in (0,1)$ and $%
\sigma ^{2}>0$. We denote by $W_{\vee }$ the pointwise maximum of the two
isotropic fractional Brownian fields, i.e. 
\begin{equation*}
W_{\vee }(x)=W^{(1)}(x)\vee W^{(2)}(x),\quad x\in \mathbf{R}^{2}.
\end{equation*}

For two distinct sites $x_{1},x_{2}\in \mathbf{R}^{2}$, let%
\begin{equation*}
U_{x_{1},x_{2}}^{(W_{\vee })}=\sigma ^{-1}d_{1,2}^{-\alpha /2}\left( W_{\vee
}(x_{2})-W_{\vee }(x_{1})\right) .
\end{equation*}%
Then we define%
\begin{eqnarray*}
V_{2,N}^{(W_{\vee })} &=&\frac{1}{\sqrt{\left\vert E_{N}\right\vert }}%
\sum_{\left( x_{1},x_{2}\right) \in E_{N}}\left( (U_{x_{1},x_{2}}^{(W_{\vee
})})^{2}-1\right) \\
V_{3,N}^{(W_{\vee })} &=&\frac{1}{\sqrt{\left\vert DT_{N}\right\vert }}%
\sum_{\left( x_{1},x_{2},x_{3}\right) \in DT_{N}}\left( \left( 
\begin{array}{cc}
U_{x_{1},x_{2}}^{(W_{\vee })} & U_{x_{1},x_{3}}^{(W_{\vee })}%
\end{array}%
\right) \left( 
\begin{array}{cc}
1 & R_{x_{1},x_{2},x_{3}} \\ 
R_{x_{1},x_{2},x_{3}} & 1%
\end{array}%
\right) ^{-1}\left( 
\begin{array}{c}
U_{x_{1},x_{2}}^{(W_{\vee })} \\ 
U_{x_{1},x_{3}}^{(W_{\vee })}%
\end{array}%
\right) -2\right) ,
\end{eqnarray*}%
where $R_{x_{1},x_{2},x_{3}}$ is given in Eq. $\left( \ref{eq:coorU}\right)
. $

The main result of this section concerns the asymptotic behaviors of $%
V_{2,N}^{(W_{\vee })}$ and $V_{3,N}^{(W_{\vee })}$. To state it, let us
denote the difference between both fractional Brownian fields as $%
W^{(2\backslash 1)}\left( x\right) =W^{(2)}\left( x\right) -W^{(1)}\left(
x\right) $ for any $x\in \mathbf{R}^{2}$. Similarly to Section 5.1 in Robert
(2020), we observe that, for any real measurable function $f:\mathbf{R}%
\rightarrow \mathbf{R}$ and for any $(x_{1},x_{2})\in E_{N}$, 
\begin{multline}
\label{eq:decompos_fU}
f(U_{x_{1},x_{2}}^{(W_{\vee })}) = f(U_{x_{1},x_{2}}^{(1)})\mathbb{I}[{%
W^{(2\backslash 1)}(x_{1})<0]}+f(U_{x_{1},x_{2}}^{(2)})\mathbb{I}[{%
W^{(2\backslash 1)}(x_{1})>0]}   \\
+\Psi _{f}\left(U_{x_{1},x_{2}}^{(1)},U_{x_{1},x_{2}}^{(2)},W^{(2\backslash
1)}(x_{1})/(\sigma d_{1,2}^{\alpha /2})\right) ,  
\end{multline}%
where 
\begin{equation*}
U_{x_{1},x_{2}}^{(1)}=\frac{1}{\sigma d_{1,2}^{\alpha /2}}\left(
W^{(1)}(x_{2})-W^{(1)}(x_{1})\right) ,\quad U_{x_{1},x_{2}}^{(2)}=\frac{1}{%
\sigma d_{1,2}^{\alpha /2}}\left( W^{(2)}(x_{2})-W^{(2)}(x_{1})\right)
\end{equation*}%
and 
\begin{equation*}
\Psi _{f}\left( x,y,w\right) =(f(y+w)-f(x))\mathbb{I}\left[ x-y\leq w\leq 0%
\right] +(f(x-w)-f(y))\mathbb{I}\left[ 0\leq w\leq x-y\right] .
\end{equation*}%
In particular, taking $f(u)=H_{2}(u)=u^{2}-1$, for all $u\in \mathbf{R}$,
and $\Psi =\Psi _{H_{2}}$, the above decomposition implies that 
\begin{equation}
V_{2,N}^{(W_{\vee })}=V_{2,N}^{(1)}+V_{2,N}^{(2)}+V_{2,N}^{(2/1)},
\label{eq:decompositionV}
\end{equation}%
where 
\begin{align*}
V_{2,N}^{(1)}& =\frac{1}{\sqrt{|E_{N}|}}\sum_{(x_{1},x_{2})\in
E_{N},W^{(2\backslash 1)}\left( x_{1}\right) <0}\left(
(U_{x_{1},x_{2}}^{(1)})^{2}-1\right) \\
V_{2,N}^{(2)}& =\frac{1}{\sqrt{|E_{N}|}}\sum_{(x_{1},x_{2})\in
E_{N},W^{(2\backslash 1)}\left( x_{1}\right) >0}\left(
(U_{x_{1},x_{2}}^{(2)})^{2}-1\right) \\
V_{2,N}^{(2/1)}& =\frac{1}{\sqrt{|E_{N}|}}\sum_{(x_{1},x_{2})\in E_{N}}\Psi
(U_{x_{1},x_{2}}^{(1)},U_{x_{1},x_{2}}^{(2)},W^{(2\backslash 1)}\left(
x_{1}\right) /(\sigma d_{1,2}^{\alpha /2})).
\end{align*}%
To obtain a similar decomposition for the triples, let us denote, for $%
-1<R<1 $, by $\Omega $ the following function%
\begin{multline}
\label{eq:def_phi} 
 \Omega (u_{1},v_{1},u_{2},v_{2},w_{1},w_{2};R)   = \frac{1}{1-R^{2}}\left[ \Psi _{H_{2}}\left( u_{1},v_{1},w_{1}\right) 
+\Psi _{H_{2}}\left( u_{2},v_{2},w_{2}\right) \right]\\
\begin{split}
& -2\frac{R}{1-R^{2}}%
\Psi _{I}\left( u_{1},v_{1},w_{1}\right) \Psi _{I}\left(
u_{2},v_{2},w_{2}\right)\\
& \left.    -2\frac{R}{1-R^{2}}\left[ u_{1}\Psi _{I}\left(
u_{2},v_{2},w_{2}\right) +u_{2}\Psi _{I}\left( u_{1},v_{1},w_{1}\right) %
\right] \mathbb{I}[w_{1}{<0]}\right.   \\
&\left. -2\frac{R}{1-R^{2}}\left[ v_{1}\Psi _{I}\left(
u_{2},v_{2},w_{2}\right) +v_{2}\Psi _{I}\left( u_{1},v_{1},w_{1}\right) %
\right] \mathbb{I}[w_{1}{>0]}\right.  
\end{split}
\end{multline}%
with $I\left( u\right) =u$ for all $u\in \mathbf{R}$. Then we have (see
Section 3.3.2 in the Supplementary Material)%\ref{cor:twotrajectories_V3} 
\begin{equation}
V_{3,N}^{(W_{\vee })}=V_{3,N}^{(1)}+V_{3,N}^{(2)}+V_{3,N}^{(2/1)},
\label{eq:decompositionV3}
\end{equation}%
where 
\begin{eqnarray*}
V_{3,N}^{(1)} &=&\frac{1}{\sqrt{\left\vert DT_{N}\right\vert }}\sum 
_{\substack{ \left( x_{1},x_{2},x_{3}\right) \in DT_{N},  \\ W^{(2\backslash
1)}\left( x_{1}\right) <0}}\left( \left( 
\begin{array}{cc}
U_{x_{1},x_{2}}^{\left( 1\right) } & U_{x_{1},x_{3}}^{\left( 1\right) }%
\end{array}%
\right) \left( 
\begin{array}{cc}
1 & R_{x_{1},x_{2},x_{3}} \\ 
R_{x_{1},x_{2},x_{3}} & 1%
\end{array}%
\right) ^{-1}\left( 
\begin{array}{c}
U_{x_{1},x_{2}}^{\left( 1\right) } \\ 
U_{x_{1},x_{3}}^{\left( 1\right) }%
\end{array}%
\right) -2\right) \\
V_{3,N}^{(2)} &=&\frac{1}{\sqrt{\left\vert DT_{N}\right\vert }}\sum 
_{\substack{ \left( x_{1},x_{2},x_{3}\right) \in DT_{N},  \\ W^{(2\backslash
1)}\left( x_{1}\right) >0}}\left( \left( 
\begin{array}{cc}
U_{x_{1},x_{2}}^{\left( 2\right) } & U_{x_{1},x_{3}}^{\left( 2\right) }%
\end{array}%
\right) \left( 
\begin{array}{cc}
1 & R_{x_{1},x_{2},x_{3}} \\ 
R_{x_{1},x_{2},x_{3}} & 1%
\end{array}%
\right) ^{-1}\left( 
\begin{array}{c}
U_{x_{1},x_{2}}^{\left( 2\right) } \\ 
U_{x_{1},x_{3}}^{\left( 2\right) }%
\end{array}%
\right) -2\right) \\
V_{3,N}^{(2/1)} &=&\frac{1}{\sqrt{\left\vert DT_{N}\right\vert }}%
\sum_{\left( x_{1},x_{2},x_{3}\right) \in DT_{N}}\Omega \left(
U_{x_{1},x_{2}}^{(1)},U_{x_{1},x_{3}}^{(1)},U_{x_{1},x_{2}}^{(2)},U_{x_{1},x_{3}}^{(2)},%
\frac{W^{(2\backslash 1)}(x_{1})}{\sigma d_{1,2}^{\alpha /2}},\frac{%
W^{(2\backslash 1)}(x_{1})}{\sigma d_{1,3}^{\alpha /2}}%
;R_{x_{1},x_{2},x_{3}}\right) .
\end{eqnarray*}

An adaptation of the proof of Theorem \ref{Prop_Conv_V_n_Gaussian} shows
that, for $\alpha \in (0,1)$, as $N\rightarrow \infty $,%
\begin{equation}
V_{2,N}^{(1)}+V_{2,N}^{(2)}\overset{\mathcal{D}}{\rightarrow }\mathcal{N}%
\left( 0,\sigma _{V_{2}}^{2}\right)  \label{eq:convVni}
\end{equation}%
and%
\begin{equation}
V_{3,N}^{(1)}+V_{3,N}^{(2)}\overset{\mathcal{D}}{\rightarrow }\mathcal{N}%
\left( 0,\sigma _{V_{3}}^{2}\right) .  \label{eq:convWni}
\end{equation}%
To obtain the asymptotic behaviors of $V_{2,N}^{(W_{\vee })}$ and $%
V_{3,N}^{(W_{\vee })}$, the asymptotic behaviors of $V_{2,N}^{(2/1)}$ and $%
V_{3,N}^{(2/1)}$ are investigated. This requires to introduce the notion of
local time of $W^{(2\backslash 1)}$.

\paragraph{The local time of $W^{(2\backslash 1)}$.}

Let $\nu ^{(2\backslash 1)}$ be the occupation measure of ${%
W^{(2\backslash 1)}}$ over $\mathbf{C}$ defined by 
\begin{equation*}
\nu ^{(2\backslash 1)}\left( A\right) =\int_{\mathbf{C}}\mathbb{I}\left[ {%
W^{(2\backslash 1)}\left( x\right) \in A}\right] \mathrm{d}x,
\end{equation*}%
for any Borel measurable set $A\subset \mathbf{R}$. Observe that, for any $%
s,t\in \lbrack 0,1]^{2}$, 
\begin{equation*}
\Delta (s,t):=\mathbb{E}\left[ {(W^{(2\backslash 1)}\left( s\right)
-W^{(2\backslash 1)}\left( t\right) )^{2}}\right] =2\sigma ^{2}\left\Vert
s-t\right\Vert ^{\alpha }.
\end{equation*}%
Because $\int_{\mathbf{C}}(\Delta (s,t))^{-1/2}\mathrm{d}s$ is finite for
all $t\in \mathbf{C}$, it follows from Section 22 in Geman and Horowitz
(1980) that the occupation measure $\nu ^{(2\backslash 1)}$ admits a
Lebesgue density, referred to as the local time, that we denote by 
\begin{equation*}
L_{W^{(2\backslash 1)}}\left( \ell \right) :=\frac{d\nu ^{(2\backslash 1)}}{%
d\ell }\left( \ell \right) .
\end{equation*}%
An immediate consequence of the existence of the local time is the
occupation time formula, which states that 
\begin{equation*}
\int_{\mathbf{C}}g(W^{(2\backslash 1)}\left( x\right) )\mathrm{d}x=\int_{%
\mathbf{R}}g\left( \ell \right) L_{W^{(2\backslash 1)}}\left( \ell \right) 
\mathrm{d}\ell
\end{equation*}%
for any Borel function $g$ on $\mathbf{R}$. Adapting the proof of Lemma 1.1
in Jaramillo et al. (2021), we can easily show that, for any $\ell \in 
\mathbf{R}$, 
\begin{equation*}
L_{W^{(2\backslash 1)}}\left( \ell \right) =\lim_{\varepsilon \rightarrow
0}\int_{\mathbf{C}}\frac{1}{\sqrt{2\pi \varepsilon }}\exp \left( -\frac{1}{%
2\varepsilon }\left( W^{(2\backslash 1)}\left( x\right) -\ell \right)
^{2}\right) \mathrm{d}x
\end{equation*}%
or 
\begin{equation}
L_{W^{(2\backslash 1)}}\left( \ell \right) =\frac{1}{2\pi }%
\lim_{M\rightarrow \infty }\int_{-[M,M]}\int_{\mathbf{R}}e^{\mathrm{i}\xi
(W^{(2\backslash 1)}\left( x\right) -\ell )}\mathrm{d}x\mathrm{d}\xi ,
\label{eq:localtimeintegral}
\end{equation}%
where the limits hold in $L^{2}$.

\paragraph{The asymptotic behaviors of $V_{2,N}^{(W_{\vee })}$ and $%
V_{3,N}^{(W_{\vee })}$.}

Let $F_{2}$ be the function defined, for any $z\in \mathbf{R}$, by 
\begin{equation*}
F_{2}(z)=\int_{\mathbf{R}^{2}\times \mathbf{R}_{+}}\Psi
_{H_{2}}(x,y,z/d^{\alpha /2})\frac{1}{2\pi }e^{-(x^{2}+y^{2})/2}f_{D}\left(
d\right) \mathrm{d}x\mathrm{d}y\mathrm{d}d,
\end{equation*}%
where $f_{D}$ is the density function of the length of the typical edge
defined in Eq. $\left( \ref{eq:typicallength}\right) $. Let us also define%
\begin{multline*}
F_{3}(z) = \int_{\mathbf{R}^{4}\times (\mathbf{R}_{+})^{3}}\Omega
(x_{1},y_{1},x_{2},y_{2},z/d_{1}^{\alpha /2},z/d_{3}^{\alpha
/2};R(d_{1},d_{2},d_{3})) \\
\begin{split}
& \times \varphi _{2}\left( x_{1},y_{1};R(d_{1},d_{2},d_{3})\right) \varphi
_{2}\left( x_{2},y_{2};R(d_{1},d_{2},d_{3})\right) \\
&\times f_{D_{1},D_{2},D_{3}}\left( d_{1},d_{2},d_{3}\right) \mathrm{d}x_{1}%
\mathrm{d}y_{1}\mathrm{d}x_{2}\mathrm{d}y_{2}\mathrm{d}d_{1}\mathrm{d}d_{2}%
\mathrm{d}d_{3},
\end{split}
\end{multline*}%
where%
\begin{eqnarray*}
\varphi _{2}\left( x,y;R\right) &=&\frac{1}{2\pi }\frac{1}{\left(
1-R^{2}\right) }\exp \left( -\frac{1}{2}\left( 
\begin{array}{cc}
x & y%
\end{array}%
\right) \left( 
\begin{array}{cc}
1 & R \\ 
R & 1%
\end{array}%
\right) ^{-1}\left( 
\begin{array}{c}
x \\ 
y%
\end{array}%
\right) \right) , \\
R(d_{1},d_{2},d_{3}) &=&\frac{d_{1}^{\alpha }+d_{3}^{\alpha }-d_{2}^{\alpha }%
}{2\left( d_{1}d_{3}\right) ^{\alpha /2}},
\end{eqnarray*}%
and where $f_{D_{1},D_{2},D_{3}}$ is the density function of the edge
lengths of the typical Delaunay triangle $\mathcal{C}$.

Moreover let%
\begin{equation*}
c_{V_{2}}=\int_{\mathbf{R}}F_{2}(z)\mathrm{d}z\quad \text{and}\quad
c_{V_{3}}=\int_{\mathbf{R}}F_{3}(z)\mathrm{d}z.
\end{equation*}%
The following proposition provides the asymptotic behaviors of $%
V_{2,N}^{(2/1)}$ and $V_{3,N}^{(2/1)}$.

\begin{proposition}
\label{prop:twotrajectories} Let $\alpha \in (0,1)$. Then, as $N\rightarrow
\infty $, 
\begin{eqnarray*}
&&\frac{\sqrt{3}}{3}N^{-(2-\alpha )/4}V_{2,N}^{(2/1)}\overset{\mathbb{P}}{%
\rightarrow }c_{V_{2}}L_{W^{(2\backslash 1)}}(0) \\
&&\frac{\sqrt{2}}{2}N^{-(2-\alpha )/4}V_{3,N}^{(2/1)}\overset{\mathbb{P}}{%
\rightarrow }c_{V_{3}}L_{W^{(2\backslash 1)}}(0).
\end{eqnarray*}
\end{proposition}

Note that the factors $\sqrt{3}/3$ and $\sqrt{2}/2$ come from the facts that 
$|E_{N}|/N\overset{a.s.}{\rightarrow }3$ and $|DT_{N}|/N\overset{a.s.}{%
\rightarrow }2$ as $N\rightarrow \infty $, respectively. As a consequence of
the above proposition, we obtain the following result.

\begin{theorem}
\label{cor:twotrajectories} Let $\alpha \in (0,1)$. Then, as $N\rightarrow
\infty $, 
\begin{eqnarray*}
&&\frac{\sqrt{3}}{3}N^{-(2-\alpha )/4}V_{2,N}^{(W_{\vee })}\overset{\mathbb{P%
}}{\rightarrow }c_{V_{2}}L_{W^{(2\backslash 1)}}(0) \\
&&\frac{\sqrt{2}}{2}N^{-(2-\alpha )/4}V_{3,N}^{(W_{\vee })}\overset{\mathbb{P%
}}{\rightarrow }c_{V_{3}}L_{W^{(2\backslash 1)}}(0).
\end{eqnarray*}
\end{theorem}

An important observation is that the rates of convergence of $%
V_{2,N}^{(W_{\vee })}$ and $V_{3,N}^{(W_{\vee })}$ differ from those of $%
V_{2,N}^{(W)}$ and $V_{3,N}^{(W)}$. The sums of square increments in $%
V_{2,N}^{(2/1)}$ and $V_{3,N}^{(2/1)}$ are actually the dominant terms.
These increments depend on both isotropic fractional Brownian fields and
they reveal the local time of $W^{(2\backslash 1)}$\ at level $0$ in the
limits. It is also noteworthy that the convergence is now in probability.

\subsection{Asymptotic distributions of squared increment sums for the
max-stable Brown-Resnick random field}

Let $\left( \eta \left( x\right) \right) _{x\in \mathbf{R}^{2}}$ be a
max-stable Brown-Resnick random field such that $\eta (x)=\bigvee_{i\geq
1}U_{i}Y_{i}(x)$ for any $x\in \mathbf{R}^{2}$, where $(U_{i})_{i\geq 1}$ is
a decreasing enumeration of the points of a Poisson point process on $%
(0,+\infty )$ with intensity measure $u^{-2}\mathrm{d}u$, and $%
(Y_{i})_{i\geq 1}$ are i.i.d. copies of \textbf{\ }%
\begin{equation*}
Y\left( x\right) =\exp \left( W\left( x\right) -\gamma \left( x\right)
\right) ,\qquad x\in \mathbf{R}^{2},
\end{equation*}%
where $\left( W\left( x\right) \right) _{x\in \mathbf{R}^{2}}$ is an
isotropic fractional Brownian field satisfying $W\left( 0\right) =0$ a.s.
and $\gamma \left( x\right) =\text{var}\left( W\left( x\right) \right) /2=\sigma
^{2}\left\Vert x\right\Vert ^{\alpha }/2$ for some $\alpha \in (0,1)$ and $%
\sigma ^{2}>0$.

Let us define, for $k\neq j\geq 1$, 
\begin{equation*}
Z_{k\backslash j}\left( x\right) =Z_{k}\left( x\right) -Z_{j}\left( x\right)
,\qquad x\in \mathbf{R}^{2},
\end{equation*}%
where%
\begin{equation*}
Z_{i}\left( x\right) =\log U_{i}+\log Y_{i}(x),\qquad x\in \mathbf{R}^{2}.
\end{equation*}%
In the same spirit as Dombry and Kabluchko (2018), we build a random
tessellation of $\mathbf{C}$, $\left( \mathbf{C}_{k,j}\right) _{k\neq j\geq
1}$ where%
\begin{equation}
\mathbf{C}_{k,j}=\left\{ x\in \mathbf{C:}Z_{k}\left( x\right) \bigwedge
Z_{j}\left( x\right) >\bigvee_{i\neq j,k}Z_{i}\left( x\right) \right\} .
\label{eq:def_Ckj}
\end{equation}%
If $\mathbf{C}_{k,j}\neq \varnothing $, we define for any Borel subset $A$
of $\mathbf{R}$ the occupation measure of $Z_{k\backslash j}$ over $%
\mathbf{C}_{k,j}$ by 
\begin{equation*}
\nu ^{(k\backslash j)}\left( A\right) =\int_{\mathbf{C}_{k,j}}\mathbb{I}%
\left[ Z_{k\backslash j}\left( x\right) {\in A}\right] \mathrm{d}x.
\end{equation*}%
The associated local time at level $0$ is given by $L_{Z_{k\backslash
j}}\left( 0\right) :=\frac{d\nu ^{(k\backslash j)}}{d\ell }\left( 0\right) $%
. If $\mathbf{C}_{k,j}=\varnothing $, we let $L_{Z_{k\backslash j}}\left(
0\right) :=0$.

Let $U_{x_{1},x_{2}}^{(\eta )}$ be the (normalized) increment of $\log
\left( \eta \right) $ defined as%
\begin{equation*}
U_{x_{1},x_{2}}^{(\eta )}=\frac{1}{\sigma \left\Vert x_{2}-x_{1}\right\Vert
^{\alpha /2}}\log \left( \frac{\eta (x_{2})}{\eta (x_{1})}\right) .
\end{equation*}%
The square increment sums are given respectively by%
\begin{eqnarray*}
V_{2,N}^{(\eta )} &=&\frac{1}{\sqrt{\left\vert E_{N}\right\vert }}%
\sum_{\left( x_{1},x_{2}\right) \in E_{N}}\left( (U_{x_{1},x_{2}}^{(\eta
)})^{2}-1\right) \\
V_{3,N}^{(\eta )} &=&\frac{1}{\sqrt{\left\vert DT_{N}\right\vert }}%
\sum_{\left( x_{1},x_{2},x_{3}\right) \in DT_{N}}\left( \left( 
\begin{array}{cc}
U_{x_{1},x_{2}}^{(\eta )} & U_{x_{1},x_{3}}^{(\eta )}%
\end{array}%
\right) \left( 
\begin{array}{cc}
1 & R_{x_{1},x_{2},x_{3}} \\ 
R_{x_{1},x_{2},x_{3}} & 1%
\end{array}%
\right) ^{-1}\left( 
\begin{array}{c}
U_{x_{1},x_{2}}^{(\eta )} \\ 
U_{x_{1},x_{3}}^{(\eta )}%
\end{array}%
\right) -2\right) ,
\end{eqnarray*}%
where $R_{x_{1},x_{2},x_{3}}$ is given in Eq. $\left( \ref{eq:coorU}\right) $%
.

\begin{theorem}
\label{prop:BRtrajectories}Let $\alpha \in (0,1)$. Then, as $N\rightarrow
\infty $, 
\begin{eqnarray*}
&&\frac{\sqrt{3}}{3}N^{-(2-\alpha )/4}V_{2,N}^{(\eta )}\overset{\mathbb{P}}{%
\rightarrow }c_{V_{2}}\sum_{j\geq 1}\sum_{k>j}L_{Z_{k\backslash j}}\left(
0\right) \\
&&\frac{\sqrt{2}}{2}N^{-(2-\alpha )/4}V_{3,N}^{(\eta )}\overset{\mathbb{P}}{%
\rightarrow }c_{V_{3}}\sum_{j\geq 1}\sum_{k>j}L_{Z_{k\backslash j}}\left(
0\right) .
\end{eqnarray*}
\end{theorem}

The results in Theorem \ref{prop:BRtrajectories} are quite similar with
those in Theorem \ref{cor:twotrajectories}. It can be noted that there is an
a.s. finite number of local times $L_{Z_{k\backslash j}}\left( 0\right) $, $%
j\geq 1$ and $k>j$, which are positive. This is related to the fact that
there is an a.s. finite number of non-empty cells of the canonical
tessellation in $\mathbf{C}$.

Using the Slivnyak-Mecke formula (see e.g. Theorem 3.2.5 in Schneider-Weil
(2008)) and the same arguments as in the proof of Proposition 3 in Robert
(2020), we can state that%
\begin{equation*}
\lim_{N\rightarrow \infty }N^{\alpha /4}\mathbb{E}\left[ \frac{1}{N}%
\sum_{\left( x_{1},x_{2}\right) \in E_{N}}\left( (U_{x_{1},x_{2}}^{(\eta
)})^{2}-1\right) \right] =4\sigma \mathbb{E}\left[ D^{\alpha /2}\right] \psi
\end{equation*}%
with%
\begin{equation*}
\psi =\int_{0}^{\infty }u\varphi (u)\left[ 1/2-\bar{\Phi}\left( u\right) -u%
\bar{\Phi}\left( u\right) \Phi \left( u\right) /\varphi (u)\right] du\simeq
-0.094.
\end{equation*}%
As a consequence we deduce that $c_{V_{2}}$ is negative.

\subsection{Asymptotic properties of the MCLEs}

We are now able to present the asymptotic properties of $\hat{\sigma}%
_{j,N}^{2}$ and $\hat{\alpha}_{j,N}$ for $j=2,3$. Let us recall that the
sums of the contributions of the observations to the composite likelihood
are proportional to the square increment statistics (see Propositions \ref%
{Prop_pairwise_pdf} and \ref{Prop_triplewise_pdf}). Moreover the asymptotic
behaviors of these statistics are characterized in Theorem \ref%
{prop:BRtrajectories}.

\begin{theorem}
\label{Prop:Asym_Prop_CL_Est}Assume that $\sigma _{0}$ belongs to the
interior of a compact set of $\mathbf{R}_{+}$, and that $\alpha _{0}$
belongs to the interior of a compact set of $(0,1)$. Then, as $N\rightarrow
\infty $,%
\begin{eqnarray*}
&&\frac{\sqrt{3}}{3}\sqrt{|E_{N}|}N^{-(2-\alpha _{0})/4}\left( \hat{\sigma}%
_{2,N}^{2}-\sigma _{0}^{2}\right) \overset{\mathbb{P}}{\rightarrow }%
c_{V_{2}}\sigma _{0}^{2}\sum_{j\geq 1}\sum_{k>j}L_{Z_{k\backslash j}}\left(
0\right) \\
&&\frac{\sqrt{3}}{6}\sqrt{|E_{N}|}N^{-(2-\alpha _{0})/4}\log (N)\left( \hat{%
\alpha}_{2,N}-\alpha _{0}\right) \overset{\mathbb{P}}{\rightarrow }%
-c_{V_{2}}\sum_{j\geq 1}\sum_{k>j}L_{Z_{k\backslash j}}\left( 0\right)
\end{eqnarray*}%
and%
\begin{eqnarray*}
&&\frac{\sqrt{2}}{2}\sqrt{|E_{N}|}N^{-(2-\alpha _{0})/4}\left( \hat{\sigma}%
_{3,N}^{2}-\sigma _{0}^{2}\right) \overset{\mathbb{P}}{\rightarrow }%
c_{V_{3}}\sigma _{0}^{2}\sum_{j\geq 1}\sum_{k>j}L_{Z_{k\backslash j}}\left(
0\right) \\
&&\frac{\sqrt{2}}{4}\sqrt{|E_{N}|}N^{-(2-\alpha _{0})/4}\log (N)\left( \hat{%
\alpha}_{3,N}-\alpha _{0}\right) \overset{\mathbb{P}}{\rightarrow }%
-c_{V_{3}}\sum_{j\geq 1}\sum_{k>j}L_{Z_{k\backslash j}}\left( 0\right) .
\end{eqnarray*}
\end{theorem}

Several important points have to be highlighted. First the MCLEs of $\sigma
_{0}^{2}$ and $\alpha _{0}$ (when the other parameter is known) are
consistent in our infill asymptotic setup. They have rates of convergence
proportional to $N^{\alpha _{0}/4}$ for $\hat{\sigma}_{2,N}^{2}$ and $\log
\left( N\right) N^{\alpha _{0}/4}$ for $\hat{\alpha}_{2,N}$ that differ from
the expected rates of convergence $N^{1/2}$ and $\log \left( N\right)
N^{1/2} $ as in Zhu and Stein (2002) for the isotropic fractional Brownian
field. Second the type of convergence is in probability. The random variables appearing in the limits in Theorem \ref{Prop:Asym_Prop_CL_Est} are proportional to a sum of local times. However these local times have unknown
distributions and they cannot be estimated from the data since the
underlying random fields $(Y_{i})_{i\geq 1}$ and the point process $\left(
U_{i}\right) _{i\geq 1}$ are not observed.
In particular, if the spatial data are only
observed for a single date, the Gaussian approximation for the MCLEs given
in Padoan et al. (2010) (when several independent replications over time of
the spatial data are available) should not be used.

The problem of joint parameter estimation of $\left( \sigma _{0}^{2},\alpha
_{0}\right) $ is left for future work, but it is expected that the
respective rates of convergence will be modified into $N^{\alpha
_{0}/4}/\log \left( N\right) $ and $N^{\alpha _{0}/4}$ as suggested by
Brouste and Fukasawa (2018) in the case of a fractional Gaussian process ($%
d=1$) observed on a regular grid.

\end{document}